\documentclass[envcountsect]{svjour3} 
\smartqed 
\usepackage{graphicx}

\usepackage{latexsym,amsmath,amsfonts,amscd,amssymb,verbatim,hyperref}
\usepackage{color}
\begin{document}

\title{Two Optimal Value Functions in Parametric Conic Linear Programming\footnote{Dedicated to Professor Franco Giannessi on the occasion of his 85th birthday.}}


\author{Nguyen Ngoc Luan \and Do Sang Kim \and Nguyen Dong Yen}

\institute{Nguyen Ngoc Luan \at Department of Mathematics and Informatics, Hanoi National University of Education\\ 
136 Xuan Thuy, Hanoi, Vietnam\\
luannn@hnue.edu.vn  
           \and
 Do Sang Kim  \at Department of Applied Mathematics, Pukyong National University\\
 Busan, Korea\\
 dskim@pknu.ac.kr
         \and
Nguyen Dong Yen \at
Institute of Mathematics, Vietnam Academy of	Science and Technology\\
18 Hoang Quoc Viet, Hanoi 10307, Vietnam\\
ndyen@math.ac.vn
}

\date{Received: date / Accepted: date}

\titlerunning{Parametric Conic Linear Programming}

\authorrunning{Nguyen Ngoc Luan, Do Sang Kim, Nguyen Dong Yen}

\maketitle

\begin{abstract}
We consider the conic linear program given by a closed convex cone in an Euclidean space and a matrix, where vector on the right-hand-side of the constraint system and the vector defining the objective function are subject to change. Using the strict feasibility condition, we prove the locally Lipschitz continuity and obtain some differentiability properties of the optimal value function of the problem under right-hand-side perturbations. For the optimal value function under linear perturbations of the objective function, similar differentiability properties are obtained under the assumption saying that both primal problem and dual problem are strictly feasible.
\end{abstract}


\section{Introduction}

If the feasible region or the objective function of a mathematical programming problem depends on a parameter, then the optimal value of the problem is a function of the parameter. In general, the optimal value function is a fairly complicated function. Continuity properties (resp., differentiability properties) of the optimal value function in parametric mathematical programming are usually classified as results on {\em stability} (resp., on {\em differential stability}) of optimization problems. Many results in this direction can be found in the books \cite[Chapters~4,5]{Bonnans_Shapiro_2000}, \cite[Chapter~4]{Mordukhovich_Book_I} and \cite[Chapter~5]{Mordukhovich_Book_II}, the papers [4-13], 
the dissertation \cite{An_Thesis}, and the references therein.

In 2001, Gauvin \cite{Gauvin_2001} investigated the sensitivity of the optimal value of a parametric linear programming problem. He gave separated formulas for the increment of the optimal value with respect to perturbations of the right-hand-side vector in the constraint system, the cost vector, and the coefficients of the matrix defining the linear constraint system. In~2016, Thieu~\cite{Thieu_2016} established some lower and upper bounds for the increment of the optimal value of a parametric linear program, where both objective vector and right-hand-side vector in the constraint system are subject to perturbations. He showed that the appearance of an extra second-order term in these estimates is a must. 

  Conic linear programming is a natural extension of linear programming. In a linear program, a feasible point usually has nonnegative components, i.e., it belongs to the nonnegative orthant of an Euclidean space. In addition, the constraints are written as linear inequalities. Meanwhile, in a conic linear program, the constraint is written as an inequality with the ordering cone being a closed convex cone in an  Euclidean space. Conic linear programs can be used in many applications that cannot be simulated by linear programs (see, e.g., \cite[Section~2.2]{BN_2001} and \cite[Chapter~6]{Luenberger_Ye_2016}). This is the reason why these special convex optimization problems have attracted much attention from researchers (see  Ben-Tal and Nemirovski~\cite{BN_2001}, Bonnans and Shapiro \cite{Bonnans_Shapiro_2000}, Luenberger and Ye~\cite{Luenberger_Ye_2016}, Chuong and Jeyakumar \cite{Chuong_Jeyakumar_2017}, and the references therein).

It is of interest to know whether results similar to those of Gauvin~\cite{Gauvin_2001} and Thieu~\cite{Thieu_2016} can be obtained for conic linear programs, or not. The aim of this paper is to show that the results of~\cite{Gauvin_2001} admit certain generalizations in conic linear programming, and some differentiability properties, as well as a locally Lipschitz continuity property, can be established  for the two related optimal value functions. The obtained results are analyzed via four concrete examples which show that the optimal value functions in a conic linear program are more complicated than that of the corresponding optimal value functions in linear programming.

The remaining part of our paper has six sections. Section~2 contains some definitions, a strong duality theorem in conic linear programming, and a lemma. Section~3 is devoted to the locally Lipschitz continuity property of the optimal value function of a conic linear program under right-hand-side perturbations. Differentiability properties of this optimal value function are studied in Section~4. Some increment formulas for the function are established in Section~5. The optimal value function of a conic linear program under linear perturbations of the objective function is studied in Section~6. The final section presents some concluding remarks.

\section{Preliminaries}
\setcounter{equation}{0}

One says that a nonempty subset $K$ of the Euclidean space $\mathbb{R}^m$ is a \textit{cone} if $tK \subset K$ for all $t\geq 0$. The (positive) \textit{dual cone} $K^*$  of a cone $K \subset \mathbb{R}^m$ is given by $K^*=\{v \in \mathbb{R}^m \; : \; v^T y \geq 0 \quad \forall y \in K\}$, where $^T$ denotes the matrix transposition. Herein, vectors of  Euclidean spaces are written as rows of real numbers in the text, but they are interpreted as columns of real numbers in matrix calculations. The closed ball centered at $a\in \mathbb{R}^m$ with radius~$\varepsilon>0$ is denoted by $\bar B(a,\varepsilon)$. By $\overline{\mathbb{R}}$ we denote the set $\mathbb R\cup\{\pm \infty\}$ of extended real values.

Let $f: \mathbb{R}^k \to \overline{\mathbb{R}}$ be a proper convex function and $\bar x \in \mathbb{R}^k$ be such that $f(\bar x)$ is finite. According to Theorem~23.1 in \cite{Rockafellar_1970}, the \textit{directional derivative}
\begin{equation*}
	f'(\bar x; h):=\displaystyle\lim\limits_{t\to 0^+} \frac{f(\bar{x}+th)-f(\bar x)}{t}
\end{equation*} of $f$ at $\bar x$ w.r.t. a direction $h \in \mathbb{R}^k$, always exists (it can take values $-\infty$ or $+\infty$), and one has $f'(\bar x; h)=\inf\limits_{t >0} \dfrac{f(\bar{x}+th)-f(\bar x)}{t}.$ 
The {\it subdifferential}  of $f$ at $\bar x$ is defined by
$\partial f(\bar x)=\big\{p \in \mathbb{R}^k \; : \; p^T(x-\bar x) \leq f(x)-f(\bar x) \ \, \forall x \in \mathbb{R}^k\big\}.$

 For a convex set $C$ in an Euclidean space $\mathbb{R}^k$, the {\it normal cone} of $C$ at $\bar x\in C$ is given by $N(\bar x; C)=\left\{p \in \mathbb{R}^k \; : \; p^T(x-\bar x) \leq 0 \ \, \forall x \in C\right\}.$

Let $F: \mathbb{R}^k \rightrightarrows \mathbb{R}^{\ell}$ be a multifunction. The \textit{domain} and the  \textit{graph} 
of $F$ are given, respectively, by the formulas ${\rm{dom}}\, F=\{x \in \mathbb{R}^k \; : \; F(x) \not = \emptyset\}$ and
$${\rm{gph}}\, F=\left\{ (x,y) \in \mathbb{R}^k \times \mathbb{R}^{\ell} \; : \; y \in F(x)\right\}.$$
If ${\rm gph}\,F$ is a convex set in $\mathbb{R}^k \times \mathbb{R}^{\ell}$, then $F$ is said to be a \textit{convex multifunction}. In that case, the \textit{coderivative} of $F$ at $(\bar x, \bar y) \in {\rm{gph}}\, F$ is the multifunction $D^* F(\bar x, \bar y): \mathbb{R}^{\ell} \rightrightarrows \mathbb{R}^{k}$ with $$D^* F(\bar x, \bar y)(v):=\left\{u \in \mathbb{R}^{k} \; : \; (u, -v) \in N \big( (\bar x, \bar y); {\rm{gph}}\, F\big)\right\}$$ for all $v \in \mathbb{R}^{\ell}$ (see, e.g., \cite[Section~1.2]{Mordukhovich_Book_I} and \cite{An_Yen_2015}).

From now on, let $K$ be a closed convex cone in $\mathbb{R}^m$. For any $y^1, y^2$ from $\mathbb{R}^m$, one writes $y^1 \geq_K y^2$ if $y^1-y^2 \in K$. Similarly, one writes $y^1>_K y^2$ if $y^1-y^2 \in {\rm int}K$, where ${\rm int}K$ denotes the interior of $K$.

Given a matrix $A \in \mathbb{R}^{m \times n}$, vectors $b \in \mathbb{R}^m$ and $c \in \mathbb{R}^n$, we consider the primal {\it conic linear optimization problem}   
\begin{equation*}
	\min \{c^Tx \; : \; x \in \mathbb{R}^n, \, Ax \geq_{K}b\}. \tag{\rm P}
\end{equation*}
Following Bental and Nemirovski~\cite[p.~52]{BN_2001}, we call
\begin{equation*}
	\max \{b^Ty \; : \; y \in \mathbb{R}^m, \, A^Ty=c, \, y\geq_{K^*}0\} \tag{\rm D}
\end{equation*} the \textit{dual problem} of (\rm P). 
The feasible region, and the solution set, and the optimal value of~(\rm P) are denoted respectively by $\mathcal{F}{\rm (P)}$, $\mathcal{S}{\rm (P)}$, and $v{\rm (P)}$. The feasible region, the solution set, and the optimal value of~(\rm D) are denoted respectively by $\mathcal{F}{\rm (D)}$, $\mathcal{S}{\rm (D)}$, and  $v({\rm D})$. By definitions, one has $$v({\rm P})=\inf\{c^Tx \; : \; x \in \mathbb{R}^n, \, Ax \geq_{K}b\}$$ and $\label{v(D)} v({\rm D})=\sup\{b^Ty \; : \; y \in \mathbb{R}^m, \, A^Ty=c, \, y\geq_{K^*}0\}.$ By a standard convention, $\inf\emptyset=+\infty$ and $\sup\emptyset=-\infty$.

Thanks to the weak duality theorem  in conic linear programming \cite[Proposition~2.3.1]{BN_2001}, one has $c^T x \geq b^Ty$ for any $x \in \mathcal{F}{\rm (P)}$ and $y \in \mathcal{F}{\rm (D)}$. Hence, $v({\rm P}) \geq v({\rm D})$. By constructing suitable examples, we can show that the strict inequality $v({\rm P}) > v({\rm D})$ is possible, i.e., a conic linear program may have a duality gap. Thus, the equality $v({\rm P}) = v({\rm D})$ is guaranteed only if one imposes a certain regularity condition. In what follows, we will rely on the regularity condition called \textit{strict feasibility}, which was intensively employed by Bental and Nemirovski~\cite{BN_2001}. 

\begin{definition}\label{strict_feassibility} {\rm (See \cite[Section 2.4]{BN_2001})} {\rm If there exists a point $x^0 \in \mathbb{R}^n$ satisfying $Ax^0 >_K b$, then one says that problem {\rm (P)} is \textit{strictly feasible}. If there is some $y^0 \in \mathbb{R}^m$ satisfying $A^Ty^0=c$  and $y^0>_{K^*}0$, then the dual problem {\rm (D)} is said to be \textit{strictly feasible}.} 
\end{definition}

The main assertion of the strong duality theorem in conic linear programming can be stated as follows.

\begin{lemma}\label{Strong_duallity}{\rm (See \cite[Theorem~2.4.1 (assertions 3a and 3b)]{BN_2001})} If {\rm (P)} is strictly feasible and $v({\rm P})>-\infty$, then {\rm (D)} has a solution and $v({\rm D})=v({\rm P})$. If {\rm (D)} is strictly feasible and one has $v({\rm D})<+\infty$, then {\rm (P)} has a solution and $v({\rm D})=v({\rm P})$. 
\end{lemma}

The following analogues of the Farkas lemma \cite[p.~200]{Rockafellar_1970} will be used repeatedly in the sequel.

\begin{lemma}\label{Farkas_type} Let there be given a closed convex cone $K$ in $\mathbb{R}^m$, a matrix $A$ in $\mathbb{R}^{m \times n}$, vectors $b \in \mathbb{R}^m$ and $c \in \mathbb{R}^n$, and a real number $\alpha$. \\
	\indent {\rm (a)} Suppose that there exists a point $x^0 \in \mathbb{R}^n$ satisfying $Ax^0 >_K b$. Then, the inequality $c^T x \geq \alpha$ is a consequence of the conic-linear inequality $Ax \geq_K b$ iff  there exists a vector $y \in K^*$ satisfying $A^Ty=c$ and $b^Ty\geq \alpha$. \\
	\indent {\rm (b)} Suppose that there exists a point $y^0 \in \mathbb{R}^m$ satisfying $A^Ty^0=c$ and $y^0 >_{K^*} 0$. Then, the inequality $b^T y \leq \alpha$ is a consequence of the system $$A^Ty=c,\  y \geq_{K^*} 0$$ iff  there exists a vector $x \in \mathbb{R}^n$ satisfying $Ax \geq_K b$ and $c^Tx\leq \alpha$. 
\end{lemma}
{\bf Proof}\ To prove assertion (a), suppose that $Ax \geq_K b$ and there exists a vector $y \in K^*$ satisfying $A^Ty=c$ and $b^Ty\geq \alpha$. Then we have
\begin{equation}\label{FK_proof_1}
	\begin{aligned}
		c^Tx - b^Ty = x^Tc - b^Ty &= x^T(A^Ty) - b^Ty\\
		&=(Ax)^Ty - b^Ty\\
		&=(Ax-b)^Ty \geq 0,  
	\end{aligned}
\end{equation}
where the last inequality in~\eqref{FK_proof_1} is valid because $Ax \geq_K b$ and $y>_{K^*}0$. Combining this with the assumption $b^Ty\geq \alpha$, yields $c^Tx \geq \alpha$. Conversely, suppose that $c^Tx \geq \alpha$ for all $x$ satisfying $Ax \geq_K b$.  It is clear that {\rm (P)} is strictly feasible and $v({\rm P})\geq\alpha$. By the first assertion of Lemma~\ref{Strong_duallity}, {\rm (D)} has a solution and one has $v({\rm P})=v({\rm D})$. Let $y$ be a solution of {\rm (D)}. Then we have $y \in K^*$, $A^Ty=c$ and $b^Ty=v({\rm D})=v({\rm P})\geq \alpha$. So, assertion (a) holds true.

Now, to prove assertion (b), suppose that $A^Ty=c, \,  y \geq_{K^*} 0$ and there exists a vector $x \in \mathbb{R}^n$ satisfying $Ax \geq_K b$ and $c^Tx\leq \alpha$. Then, we have
\begin{equation}\label{FK_proof_2}
	\begin{aligned}
		b^Ty -c^Tx = b^Ty -x^Tc &= b^Ty- x^T(A^Ty) \\
		&=b^Ty-(Ax)^Ty\\
		&=-(Ax-b)^Ty \leq 0,  
	\end{aligned}
\end{equation}
where the last inequality in~\eqref{FK_proof_2} is valid because $Ax \geq_K b$ and $y\geq_{K^*}0$. Since $c^Tx \leq \alpha$, this yields $b^Ty\leq \alpha$. Conversely, suppose that $b^Ty \leq \alpha$ for all $y$ satisfying $A^Ty=c, \,  y \geq_{K^*} 0$.  So, {\rm (D)} is strictly feasible by our assumption and $v({\rm D})\leq\alpha$. By the second assertion of Lemma~\ref{Strong_duallity}, {\rm (P)} has a solution and $v({\rm P})=v({\rm D})$. Let $x$ be a solution of {\rm (P)}. Then we have $Ax \geq_K b$ and $c^Tx=v({\rm P})=v({\rm D})\leq \alpha$. This justifies the validity of (b).
$\hfill \Box$

\begin{remark} {\rm The result in Lemma~\ref{Farkas_type}(a) is a known one (see \cite[Proposition~2.4.3]{BN_2001}). The assertions of Lemma~\ref{Strong_duallity} and Lemma~~\ref{Farkas_type}(a) may be false if the assumption on the strict feasiblity of the conic-linear inequality $Ax \geq_K b$ is removed (see \cite[Example~2.4.2]{BN_2001}). Similarly, the assertion of Lemma~\ref{Farkas_type}(b) may be false if the assumption on the strict feasiblity of the conic-linear inequality $y\geq_{K^*} 0$ is removed.}
\end{remark}       

The optimal value $v({\rm P})$ depends on the parameters $A$, $b$, and~$c$. Following Gauvin~\cite{Gauvin_2001}, who studied differential stability properties of  linear programming problems, we consider the functions  $\varphi(b):=\inf\{c^Tx \; : \; x \in \mathbb{R}^n, \, Ax \geq_{K}b\}$ 
and $\psi(c):=\inf\{c^Tx \; : \; x \in \mathbb{R}^n, \, Ax \geq_{K}b\}.$
Note that $\varphi(\cdot)$ is the \textit{optimal value function of} (P) \textit{under right-hand-side perturbations of the constraint set} and $\psi(\cdot)$ is the \textit{optimal value function of} (P) \textit{under perturbations of the objective function.}

\section{Locally Lipschitz Continuity of $\varphi(\cdot)$}

According to a remark given in the paper of An and Yen~\cite[p.~113]{An_Yen_2015}, we know that~$\varphi(\cdot)$ is a convex function. The next proposition presents additional continuity properties of~$\varphi$. 

\begin{theorem}\label{Lipschitz_continuity} Suppose that {\rm (P)} is strictly feasible and $v(\rm P)>-\infty$. Then, there exists $\varepsilon>0$ such that $\varphi(b')$ is finite for every $b'\in \bar B(b,\varepsilon)$ and $\varphi$ is Lipschitz on $\bar B(b,\varepsilon)$. 
\end{theorem}
{\bf Proof}\  On one hand, since {\rm (P)} is strictly feasible, one selects a point $x^0 \in \mathbb{R}^n$ satisfying $Ax^0 >_K b$, i.e., $Ax^0 - b \in {\rm int}K$. Then there exists a convex open neighborhood of $b$ satisfying $Ax^0 - b' \in {\rm int}K$ for every $b'\in V$. So, $x^0$ is a feasible point of the problem 
\begin{equation}\label{P_b'}
	\min \{c^Tx \; : \; x \in \mathbb{R}^n, \, Ax \geq_{K}b'\} \tag{${\rm P}_{b'}$}
\end{equation}
for all $b'\in V$. Moreover $\varphi(b') \leq c^Tx^0 < +\infty$. On the other hand, applying the weak duality theorem \cite[Proposition~2.3.1]{BN_2001} for the problems \eqref{P_b'}
and
\begin{equation}\label{D_b'}
	\max \{(b')^Ty \; : \; y \in \mathbb{R}^m,\, A^Ty=c,\, y\geq_{K^*}0\} \tag{${\rm D}_{b'}$}
\end{equation}   
one has $c^Tx \geq (b')^Ty$ for all $x \in \mathbb{R}^n$ satisfying $Ax \geq_K b'$ and for every $y \in \mathbb{R}^m$ satisfying $A^Ty=c$, $y \geq_{K^*}0$. Therefore,
\begin{equation*}
	\inf \{c^Tx \; : \; x \in \mathbb{R}^n, \, Ax \geq_{K}b'\} \geq \sup \{(b')^Ty \; : \; y \in \mathbb{R}^m,\, A^Ty=c,\, y\geq_{K^*}0\}.
\end{equation*}
This means that $$
\varphi(b') \geq \sup \{(b')^Ty \; : \; y \in \mathbb{R}^m,\, A^Ty=c,\, y\geq_{K^*}0\}.$$
So, we have $\varphi(b') > -\infty$ because
\begin{equation}\label{three_equalities_1}\begin{array}{rl}
& \sup \{(b')^Ty \; : \; y \in \mathbb{R}^m,\, A^Ty=c,\, y\geq_{K^*}0\} \\ & \geq \sup \{(b')^Ty \; : \; y \in \mathcal{S}({\rm D})\} > -\infty,
\end{array}
\end{equation}
where the last inequality in~\eqref{three_equalities_1} is valid as $\mathcal{S}({\rm D})$ is nonempty (see Lemma~\ref{Strong_duallity}). It follows that $\varphi(b')$ is finite for every $b' \in V$, i.e., $V \subset {\rm int}({\rm dom}\varphi)$. Combining this with the convexity of $\varphi$,  by \cite[Theorem~24.7]{Rockafellar_1970} we can asserts that $\varphi$ is locally Lipschitz on $V$.      
$\hfill \Box$ 

\section{Differentiability Properties of $\varphi(\cdot)$}

First, let us show that the solution set of (D) possesses a remarkable stability property, provided that the assumptions of Theorem~\ref{Lipschitz_continuity} are satisfied.

\begin{proposition}\label{compact_solution} Suppose that {\rm (P)} is strictly feasible and $v(\rm P)>-\infty$. Then there exists $\varepsilon>0$ such that $\mathcal{S}({\rm D}_{b'})$ is nonempty and compact for every $b'$ in $B(b,\varepsilon)$.
\end{proposition}
{\bf Proof} Since {\rm (P)} is strictly feasible, there is $x^0 \in \mathbb{R}^n$ satisfying $Ax^0-b \in {\rm int}K$. Therefore, there is a real number $\varepsilon>0$ satisfying $Ax^0- b' \in {\rm int}K$ for all $b'\in B(b,\varepsilon)$. According to Theorem~\ref{Lipschitz_continuity}, without loss of generality we can assume that $\varphi(b')$ is finite for all $b'\in B(b,\varepsilon)$. Fix a vector $b'\in  B(b,\varepsilon)$. Then, the problem \eqref{P_b'}
is strictly feasible and $v\eqref{P_b'}>-\infty$. So, applied for the pair problems \eqref{P_b'} and \eqref{D_b'} Lemma~\ref{Strong_duallity} asserts that $\mathcal{S}\eqref{D_b'}$ is nonempty and
\begin{equation*}
	\mathcal{S}\eqref{D_b'} =\big\{y \in \mathcal{F}\eqref{D_b'} \; : \; (b')^Ty=\varphi(b')\big\}.
\end{equation*}  
Since $\mathcal{F}\eqref{D_b'}=\mathcal{F}({\rm D})$, this implies  $\mathcal{S}\eqref{D_b'}=\big\{y \in \mathcal{F}({\rm D}) \; : \; (b')^Ty=\varphi(b')\big\}.$
Clearly, $\mathcal{S}\eqref{D_b'}$ is a closed set. To prove that $\mathcal{S}\eqref{D_b'}$ is compact by contradiction, let us suppose that $\mathcal{S}\eqref{D_b'}$ is unbounded. Select a sequence $\{y^k\}$ in $\mathcal{S}\eqref{D_b'}$ satisfying the condition $\lim\limits_{k \to\infty}\|y^k\|=+\infty$. Define $\widetilde{y}^k=\|y^k\|^{-1}y^k$. Since $\|\widetilde{y}^k\|=1$, without loss of generality we can assume that the  sequence $\{\widetilde{y}^k\}$ converges to a vector $\widetilde{y}$ with $\|\widetilde{y}\|=1$.  Since the sequence $\{y^k\}$ is contained in the closed  cone~$K^*$, one has $\widetilde{y} \in K^*$. Observe that
\begin{equation*}
	\begin{aligned}
		A^T\widetilde{y}=\lim\limits_{k \to \infty}A^T\widetilde{y}^k &=\lim\limits_{k \to \infty}\dfrac{A^Ty^k}{\|y^k\|}\\
		&=\lim\limits_{k \to \infty}\dfrac{c}{\|y^k\|}=0
	\end{aligned}
\end{equation*}
and
\begin{equation*}
	\begin{aligned}
		(b')^T\widetilde{y}=\lim\limits_{k \to \infty}(b')^T\widetilde{y}^k & =\lim\limits_{k \to \infty}\dfrac{(b')^Ty^k}{\|y^k\|}\\
		&=\lim\limits_{k \to\infty}\dfrac{\varphi(b')}{\|y^k\|}=0.
	\end{aligned}
\end{equation*}
It follows that $(Ax^0-b')^T\widetilde{y}={(x^0)}^TA^T\widetilde{y}-(b')^T\widetilde{y}=0$. Meanwhile, since $Ax^0-b' \in {\rm int}K$ and $\widetilde{y} \in K^*\setminus\{0\}$, $(Ax^0-b')^T\widetilde{y} >0$. We have obtained a contradiction. Thus,  $\mathcal{S}\eqref{D_b'}$ is bounded; hence it is  compact. $\hfill \Box$

The following question arises naturally from Proposition~\ref{compact_solution}: \textit{If {\rm (D)} is strictly feasible and $v(\rm D)<+\infty$, then the solution set $\mathcal{S}({\rm P})$ is nonempty and compact, or not?} By the second assertion of Lemma~\ref{Strong_duallity}, if {\rm (D)} is strictly feasible and $v(\rm D)<+\infty$, then $\mathcal{S}({\rm P})$ is nonempty. However, $\mathcal{S}({\rm P})$ may not be a compact set. The next example is an illustration for this statement.

\begin{example} {\rm Consider the problem (P) with $n=m=2$, $K=\mathbb{R}^2_+$, $A=
		\begin{bmatrix}
			1 & 0\\
			0 & 0
		\end{bmatrix}$, $b=(1,-1)^T$, $c=(1,0)^T$. Here we have $K^*=K=\mathbb{R}^2_+$. For $y^0:=(1,1)^T$, one has $A^Ty^0=c$  and $y^0>_{K^*}0$. Thus, {\rm (D)} is strictly feasible. Meanwhile, $\mathcal{S}({\rm P})=\{1\}\times\mathbb R$ is unbounded.} 
\end{example}

The following result gives a formula for the directional derivative of $\varphi$ at~$b$ in a direction~$d$.

\begin{theorem}\label{directional_derivative} Suppose that {\rm (P)} is strictly feasible and {\rm (P)} has a solution. Then, for every direction $d \in \mathbb{R}^m$, one has 
	\begin{equation}\label{diff_varphi}
		\varphi'(b;d)=\max\limits_{y \in \mathcal{S}({\rm D})} y^Td.
	\end{equation}	 
\end{theorem}
{\bf Proof} Fix a direction $d \in \mathbb{R}^m$, consider the function $g: \mathbb{R} \to \overline{\mathbb{R}}$ given by
\begin{equation*}
	g(t)=\inf\{c^Tx \; : \; x \in \mathbb{R}^n, \, Ax \geq_{K}b+td\},
\end{equation*}
where $t \in \mathbb{R}$. Let $f(x, t)=c^Tx$ for all $(x,t)\in \mathbb{R}^n \times \mathbb{R}$. Clearly, $f$ is convex and continuous on $\mathbb{R}^n \times \mathbb{R}$. Define the multifunction $G: \mathbb{R} \rightrightarrows \mathbb{R}^n$ by setting $$G(t)=\big\{x \in \mathbb{R}^n \; : \; Ax \geq_{K}b+td \big\}.$$ Note that $G$ is a convex multifunction. Consider the parametric convex optimization problem 
$\min \{f(x,t)\; : \; x\in G(t)\},$ which depends on the parameter $t\in\mathbb R$, and observe that $g(\cdot)$ is the optimal value function of this problem. By a remark given in~\cite[p.~113]{An_Yen_2015} we know that $g(\cdot)$ is a convex function. Applying a result of An and Yen~\cite[ Theorem~4.2]{An_Yen_2015} for $\bar{t}:=0$ and for a solution $\bar{x}$ of~$({\rm P})$, we have
\begin{equation*}
	\partial g(0)=\bigcup\limits_{(x^*,t^*) \in \partial f(\bar x, \bar t)}   \big\{t^* + D^*G( \bar t, \bar x)(x^*) \big\}.
\end{equation*}
Since $\partial f(x,t)=\left\{\begin{bmatrix}
	c \\ 0
\end{bmatrix}\right\}$, this equality yields
\begin{equation}\label{subdifferential_of_g}\partial g(0)=D^*G(0, \bar x)(c).\end{equation} By the definition of coderivative for convex multifunctions, 
\begin{equation}\label{sub_diff_g}
	\begin{aligned}
		D^*G(0, \bar x)(c)&=\{t^* \in \mathbb{R} \; : \; (t^*, -c) \in N((0,\bar x), \rm{gph}\, G)\}\\
		&=\big\{t^* \in \mathbb{R} \; : \; t^* t - c^T(x - \bar x) \leq 0 \ \, \forall (t,x) \in \rm{gph}\, G\big\}\\
		&=\big\{t^* \in \mathbb{R} \; : \; c^Tx -t^* t \geq  c^T\bar x \ \, \forall (t,x) \textrm{ with } Ax \geq_K b+td \big\}\\
		&=\big\{t^* \in \mathbb{R} \; : \; c^Tx -t^* t \geq  c^T\bar x \ \, \forall (x,t) \textrm{ with } Ax -td \geq_K b \big\}.\\
	\end{aligned}
\end{equation} 
Let
\begin{equation*}
	\widetilde{c}:=\begin{bmatrix}
		c\\
		-t^*
	\end{bmatrix}, \quad \widetilde{A}:=\begin{bmatrix}
		A \quad -d
	\end{bmatrix}, \quad \widetilde{x}:=\begin{bmatrix}
		x\\
		t
	\end{bmatrix}, \quad \alpha:=c^T\bar{x}.
\end{equation*}
Clearly, the inequality $Ax-td \geq_K b$ is equivalent to $\widetilde{A}\widetilde{x} \geq_K b$. So, from \eqref{sub_diff_g} one gets
\begin{equation}\label{sub_diff_g1}
	\begin{aligned}
		D^*G(0, \bar x)(c)=\big\{t^* \in \mathbb{R} \; : \; \widetilde{c}^T \widetilde{x} \geq  \alpha \ \, \forall \widetilde{x} \in \mathbb{R}^n \times \mathbb{R} \textrm{ with } \widetilde{A}\widetilde{x} \geq_K b \big\}.\\
	\end{aligned}
\end{equation} 
Since {\rm (P)} is strictly feasible, there exists $x^0 \in \mathbb{R}^n$ satisfying $Ax^0>_K b$. Setting $\widetilde{x}^0=\begin{bmatrix}
	x^0\\
	0
\end{bmatrix}$, we have $\widetilde{A}\widetilde{x}^0 >_K b$. Then, applying the Farkas-type result in Lemma~\ref{Farkas_type}(a) to the inequality system $\widetilde{A}\widetilde{x} \geq_K b$ and the vector $\widetilde{c}$, we can assert that $\widetilde{c}^T \widetilde{x} \geq  \alpha$ for all $\widetilde{x} \in  \mathbb{R}^n \times \mathbb{R}$ such that $\widetilde{A}\widetilde{x} \geq_K b$ iff  there exists a vector $y \in K^*$ satisfying $\widetilde{A}^Ty=\widetilde{c}$ and $b^Ty \geq \alpha$. Since $$\widetilde{A}^Ty = \begin{bmatrix}
	A^T \\ -d^T
\end{bmatrix}y=\begin{bmatrix}
	A^Ty \\ -d^Ty
\end{bmatrix},$$ the equality $\widetilde{A}^Ty=\widetilde{c}$ is equivalent to the conditions $A^Ty=c$ and $d^Ty=t^*$.
Therefore, combining~\eqref{sub_diff_g1} with the description of the feasible region of the dual problem ({\rm D}), we have
\begin{equation}\label{co_diff_G}
	\begin{aligned}
		D^*G(0, \bar x)(c)&=\big\{t^* \in \mathbb{R} \; : \; \exists y \in K^* \; \textrm{s.t.} \; \widetilde{A}^Ty=\widetilde{c},\; b^Ty\geq \alpha\big\}\\
		&=\big\{t^* \in \mathbb{R} \; : \; \exists y \in K^* \; \textrm{s.t.} \; A^Ty=c,\; d^Ty=t^*,\; b^Ty\geq c^T\bar x\big\}\\
		&=\big\{d^Ty \; : \; y \in K^*,\; A^Ty=c,\; b^Ty\geq c^T\bar x \big\}\\
		&=\big\{d^Ty \; : \; y \in \mathcal{F}({\rm D}),\;  b^Ty\geq c^T\bar x\big\}\\
		&=\big\{d^Ty \; : \; y \in \mathcal{F}({\rm D}),\;  b^Ty\geq v({\rm P})\big\}.
	\end{aligned}
\end{equation} 
Since {\rm (P)} is strictly feasible and {\rm (P)} has a solution, by the strong duality theorem in Lemma~\ref{Strong_duallity}, (D) has a solution and $v({\rm D})=v({\rm P})$. Hence, by \eqref{co_diff_G},
\begin{equation}\label{co_diff_G1}
	\begin{aligned}
		D^*G(0, \bar x)(c)&=\big\{d^Ty \; : \; y \in \mathcal{F}({\rm D}),\;  b^Ty\geq v({\rm D})\big\}\\
		&=\big\{d^Ty \; : \; y \in \mathcal{F}({\rm D}),\;  b^Ty= v({\rm D})\big\}\\
		&=\big\{d^Ty \; : \; y \in \mathcal{S}({\rm D})\big\}.
	\end{aligned}
\end{equation}   
Combining \eqref{co_diff_G1} with \eqref{subdifferential_of_g} gives $\partial g(0)=\{d^Ty \; : \; y \in \mathcal{S}({\rm D})\}.$
Then, by the well-known result on the relationships between the directional derivative of a convex function \cite[Theorem~23.4]{Rockafellar_1970}, we obtain
\begin{equation}\label{three_equalities}
	g'(0;1)=\sup\limits_{s \in \partial g(0)}s=\sup\limits_{y \in \mathcal{S}({\rm D})} y^Td=\max\limits_{y \in \mathcal{S}({\rm D})}y^Td,
\end{equation} where the last equality in~\eqref{three_equalities} is valid because $\mathcal{S}({\rm D})$ is compact. Finally, using~\eqref{three_equalities} and the fact that $\varphi'(b;d)=g'(0;1)$, we get the desired formula~\eqref{diff_varphi}.~$\hfill \Box$

Based on Theorem~\ref{directional_derivative}, the next proposition provides us with a formula for the subdifferential of $\varphi$ at $b$. 
\begin{proposition}\label{subdifferential} If {\rm (P)} is strictly feasible and {\rm (P)} possesses a solution, then $\partial\varphi(b)=\mathcal{S}({\rm D})$. 
\end{proposition}
{\bf Proof}  Take a vector $p \in \partial\varphi(b)$. For any $d \in \mathbb{R}^m$ and $t>0$, from the inequality $\varphi(b+td)-\varphi(b) \geq p^T (td)$ it follows that
$\dfrac{1}{t}\left(\varphi(b+td)-\varphi(b)\right) \geq p^T d.$ Letting $t \to 0^+$, one has $\varphi'(b;d) \geq p^Td.$ Combining this with the equation~\eqref{diff_varphi} gives 
\begin{equation}\label{inequality_1}
	\max\limits_{y \in \mathcal{S}({\rm D})} y^Td \geq p^Td.
\end{equation}
To show that $p$ belongs to $\mathcal{S}({\rm D})$, suppose the contrary: $p \notin \mathcal{S}({\rm D})$. By  $\mathcal{S}({\rm D})$ is a nonempty convex set, by the strongly separation theorem (see, e.g., \cite[Corollary~11.4.2]{Rockafellar_1970}), there exists a vector $d \in \mathbb{R}^m$ such that $\max\limits_{y \in \mathcal{S}({\rm D})} y^Td < p^Td$. This contradicts with~\eqref{inequality_1}. We have thus proved that $\partial\varphi(b) \subset \mathcal{S}({\rm D})$. To obtain the opposite inclusion, take any $p \in \mathcal{S}({\rm D})$. For all $d \in \mathbb{R}^m$, by \eqref{diff_varphi}, 
\begin{equation}\label{inequality_2}
	\varphi'(b;d)=\max\limits_{y \in \mathcal{S}({\rm D})} y^Td  \geq p^T d.
\end{equation}  
Besides, by \cite[Theorem~23.1]{Rockafellar_1970} one gets
\begin{equation}\label{inequality_3}
	\varphi(b+d) - \varphi(b) \geq \varphi'(b;d).
\end{equation} 
Combining \eqref{inequality_2} with \eqref{inequality_3} we have $\varphi(b+d) - \varphi(b) \geq p^Td$ for all $p \in \mathbb{R}^m$. So, $p \in \partial\varphi(b)$. The proof is complete. $\hfill \Box$

\section{Increment Estimates for $\varphi(\cdot)$}

Roughly speaking, the following statement gives an analogue of \cite[Theorem~1]{Gauvin_2001} for conic linear programs.

\begin{theorem}\label{Thm1} Suppose that {\rm (P)} is strictly feasible and $v{\rm (P)}>-\infty$. Then, for every $d \in \mathbb{R}^n$ and $t>0$, one has
	\begin{equation*}
		\varphi(b+td) \geq \varphi(b)+ t \sup \{d^Ty \; : \; y \in \mathbb{R}^m,\ \, A^Ty=c,\  y\geq_{K^*}0,\  b^Ty=\varphi(b)\}.
	\end{equation*}
	In addition, for every $d \in \mathbb{R}^n$, there exists $\tau >0$ such that 
	\begin{equation}\label{upper_ess}
		\varphi(b+td) \leq \varphi(b)+ t \sup \{d^Ty \; : \; y \in \mathcal{F}{\rm (D)}\}
	\end{equation}	     
	for all $t \in (0,\tau]$. 
\end{theorem}
{\bf Proof} By Lemma~\ref{Strong_duallity}, (D) has a solution and $v({\rm D})=v({\rm P})$; hence $v({\rm D})=\varphi(b)$. On one hand, applying the weak duality theorem \cite[Proposition~2.3.1]{BN_2001} for the problems 
\begin{equation}\label{P_t}
	\min \{c^Tx \; : \; x \in \mathbb{R}^n, \, Ax \geq_{K}b+td\} \tag{${\rm P}_t$}
\end{equation}
and
\begin{equation}\label{D_t}
	\max \{(b+td)^Ty \; : \; y \in \mathbb{R}^m,\, A^Ty=c,\, y\geq_{K^*}0\} \tag{${\rm D}_t$}
\end{equation}   
one has $c^Tx \geq (b+td)^Ty$ for all $x \in \mathbb{R}^n$ satisfying $Ax \geq_K b+td$ and for every $y \in \mathbb{R}^m$ satisfying $A^Ty=c$, $y \geq_{K^*}0$. Consequently, 
\begin{equation}\label{ineq_1}\begin{array}{rl}
	& \inf\big\{c^Tx \; : \; x \in \mathbb{R}^n, Ax \geq_K b+td\big\}\\ &\geq \sup\big\{(b+td)^Ty \; : \; y \in \mathbb{R}^m, A^Ty=c, y \geq_{K^*}0\big\}.
\end{array}
\end{equation}     
On the other hand, it is clear that
\begin{multline}\label{ineq_2}
	\sup\{(b+td)^Ty \; : \; y \in \mathbb{R}^m, A^Ty=c, y \geq_{K^*}0\} \\
	\geq  \sup\{(b+td)^Ty \; : \; y \in \mathbb{R}^m, A^Ty=c, y \geq_{K^*}0,b^Ty=\varphi(b)\}.
\end{multline}
Combining \eqref{ineq_1} and \eqref{ineq_2}, one gets
\begin{multline*}
	\inf\{c^Tx \; : \; x \in \mathbb{R}^n, Ax \geq_K b+td\} \\
	\geq  \sup\{(b+td)^Ty \; : \; y \in \mathbb{R}^m, A^Ty=c, y \geq_{K^*}0,b^Ty=\varphi(b)\}.
\end{multline*}
It follows that
\begin{equation}\label{lower_ineq}
	\begin{array}{rl}
		\varphi(b+td) &\geq  \sup\{(b+td)^Ty \; : \; y \in \mathbb{R}^m, A^Ty=c, y \geq_{K^*}0,b^Ty=\varphi(b)\}\\
		&=\sup\{b^Ty+td^Ty\; : \; y \in \mathbb{R}^m, A^Ty=c, y \geq_{K^*}0,b^Ty=\varphi(b)\}\\
		&=\varphi(b)+t\sup\{d^Ty\; : \; y \in \mathbb{R}^m, A^Ty=c, y \geq_{K^*}0,b^Ty=\varphi(b)\}.
	\end{array}
\end{equation}
Combining \eqref{lower_ineq} with the fact that $\varphi(b)=v({\rm P}) > -\infty$, one has $\varphi(b+td) >  -\infty$.

Since (P) is strictly feasible, there is a point $x^0 \in \mathbb{R}^n$ satisfying $Ax^0 >_K b$. Then one can find $\tau >0$ such that $Ax^0  >_K b+td$ for all $t \in [0,\tau]$. Hence, (${\rm P}_t$) is strictly feasible. Applying Lemma~\ref{Strong_duallity} for a pair problems (${\rm P}_t$) and (${\rm D}_t$), where $t \in [0,\tau]$, we have
\begin{equation*}
	\begin{aligned}
		\varphi(b+td) &= \sup\{(b+td)^Ty \; : \; y \in \mathbb{R}^m, A^Ty=c, y \geq_{K^*}0\}\\
		& \leq \sup\{b^Ty \; : \; y \in \mathbb{R}^m, A^Ty=c, y \geq_{K^*}0\}\\
		&\hspace*{0.5cm} + \sup\{td^Ty \; : \; y \in \mathbb{R}^m, A^Ty=c, y \geq_{K^*}0\}\\
		&=\varphi(b)+ t\sup\{d^Ty \; : \; y \in \mathbb{R}^m, A^Ty=c, y \geq_{K^*}0\}.
	\end{aligned}
\end{equation*}     
Thus, the inequality \eqref{upper_ess} has been proved for all $t \in [0,\tau]$. $\hfill \Box$

We now consider the case where $K$ is a polyhedral convex cone.

\begin{theorem}\label{PCC_Case} Suppose that $K$ is a polyhedral convex cone and $\varphi(b)$ is finite. Then for every $d \in \mathbb{R}^n$, there exists $\tau >0$ such that 
	\begin{equation}\label{eq_ess}
		\varphi(b+td) = \varphi(b)+ t \sup \{d^Ty \; : \; y \in \mathbb{R}^m, \, A^Ty=c, \, y\geq_{K^*}0, b^Ty=\varphi(b)\}
	\end{equation}
	for all $t \in ]0,\tau]$. 
\end{theorem}
{\bf Proof} Since $K$ is a polyhedral convex cone in $\mathbb{R}^m$, there exists a matrix $B \in \mathbb{R}^{k \times m}$ satisfying 
$K=\{y \in \mathbb{R}^m \; : \; By \geq_{\mathbb{R}^k_+} 0\}.$
Then, the problem {\rm (P)} can be written in the form
\begin{equation*}
	\min \{c^Tx \; : \; x \in \mathbb{R}^n, \, BAx \geq_{\mathbb{R}^k_+} Bb\}. \tag{\rm P'}
\end{equation*}
The dual problem of ({\rm P'}) is the problem 
\begin{equation*}
	\max \{(Bb)^Tz \; : \; z \in \mathbb{R}^k, \, (BA)^Tz=c, z \geq_{\mathbb{R}^k_+} 0\}. \tag{\rm D'}
\end{equation*}
Since $\varphi(b)$ is finite, the ({\rm P'}) has a solution. Moreover, both of two problems ({\rm P'}) and~({\rm D'}) have solutions and the optimal values are equal. Then, for every $d \in \mathbb{R}^m$, by using Theorem~1 in~\cite{Gauvin_2001} for the problem ({\rm P'}), one gets
\begin{equation*}
	\begin{aligned}
		\varphi(b+td)&=\inf\{c^Tx \; : \; x \in \mathbb{R}^n, \,  Ax \geq_K b+td \}\\
		&=\inf\{c^Tx \; : \; x \in \mathbb{R}^n, \, BAx \geq_{\mathbb{R}^k_+} Bb+tBd\}\\
		&=\varphi(b)+t \sup\{(Bd)^Tz \; : \; z \in \mathbb{R}^k_+, (BA)^Tz=c, (Bb)^Tz=\varphi(b)\}\\
		&=\varphi(b)+t \sup\{d^T(B^Tz) \; : \; z \in \mathbb{R}^k_+, A^T(B^Tz)=c, b^T(B^Tz)=\varphi(b)\}\\
	\end{aligned}
\end{equation*}
for any $t>0$ sufficiently small. Denote $B^Tz=y$, we have
\begin{equation}\label{ess_eq_1}
	\varphi(b+td)=\varphi(b)+t \sup\{d^Ty \; : \; y \in B^T(\mathbb{R}^k_+), A^Ty=c, b^Ty=\varphi(b)\}.\\
\end{equation}
In the other hand, since $K=\{y \in \mathbb{R}^m \; : \; \langle B_i^T, y \rangle \geq 0, \ i=1,2,\dots,k  \}$
where $B_i$ is the $i$-row of matrix $B$, one has 
$$K^*=\left\{\sum\limits_{i=1}^k \lambda_i B^T_i \; : \; \lambda_i \geq 0,\, i=1,2,\dots,k \right \}$$
by Proposition~2.42 in \cite{Bonnans_Shapiro_2000}; hence $K^*=B^T(\mathbb{R}^k_+)$. Combining the later and the equation~\eqref{ess_eq_1}, we can assert that \eqref{eq_ess} holds.     
$\hfill \Box$

\begin{remark}{\rm Theorem~\ref{PCC_Case} is a generalization of Theorem~1 in \cite{Gauvin_2001}, where the case $K=\mathbb{R}^m_+$ was treated.}   
\end{remark}

The following two examples are designed as illustrations for Theorem~\ref{Thm1}.

\begin{example}\label{Exam4}{\rm Consider the conic problem (P) with $x=(x_1, x_2)^T \in \mathbb{R}^2$, $K$ is the 3D ice cream cone in $\mathbb{R}^3$, i.e., 
		$K=\left\{(y_1, y_2, y_3)^T \in\mathbb{R}^3\; : \; y_3 \geq \sqrt{y_1^2+y_2^2}\right\}$, $A=
		\begin{bmatrix}
			0 & 0\\
			0 & 1\\
			1 & 0
		\end{bmatrix}$, $b=(-1,0,0)^T$, $c=(1,0)^T$. We have	
		\begin{equation}\label{exam_feasible_region}\mathcal{F}({\rm P})=\Big\{x \in \mathbb{R}^2 \; : \; x_1 \geq \sqrt{1+x_2^2}\Big\}.\end{equation} Clearly, $\varphi(b)=1$ and $\mathcal{S}({\rm P})=\big\{(1,0)^T\big\}.$ Note that $A^Ty=c$ iff  $y_2=0$ and $y_3=1$. In addition $K^*=K$. It follows that $$\mathcal{F}({\rm D})=\{y=(y_1,0,1)^T \in \mathbb{R}^3 \; : \;  -1 \leq y_1 \leq 1\}.$$
		Obviously, $\mathcal{S}({\rm D})=\{(-1,0,1)^T\}$. Given $d \in \mathbb{R}^3$ with $d_1>0$, and for every $t >0$, one has $Ax \geq_K b+td$ iff  $x_1 \geq \sqrt{(1-td_1)^2+(x_2-td_2)^2}+td_3$. On one hand, for~$t$ small enough, 
		\begin{equation*}
			\varphi(b+td)=\sqrt{(1-td_1)^2}+td_3=|1-td_1|+td_3=1+t(d_3-d_1).
		\end{equation*}
		On the other hand, $
		\varphi(b)+ t \sup \{d^Ty \; : \; y \in \mathcal{S}({\rm D})\}=1+t(d_3-d_1)$
		and
		\begin{equation*}
			\begin{aligned}
				\varphi(b)+ t \sup \{d^Ty \; : \; y \in \mathcal{F}({\rm D})\}&=1+t \sup \{d_1y_1+d_3 \; : \;  -1 \leq y_1 \leq 1\}\\
				&=1+t(|d_1|+d_3).
			\end{aligned}
		\end{equation*}  	 
}\end{example}

\begin{example}{\rm Consider the conic problem (P) with $x\in \mathbb{R}$, $K$ is the 3D ice cream cone in $\mathbb{R}^3$, $A=
		\begin{bmatrix}
			1 \\
			0\\
			0
		\end{bmatrix}$, $b=(0,0,-1)^T$, $c=1$, and $d=(0,1,-1)^T$. Clearly, the condition $Ax \geq_K b+td$ is equivalent to $1+t \geq \sqrt{x^2+t^2}$. For $t>0$ is small enough, the latter can be rewritten as $x^2 \leq 2t+1$. It follows that
		$\mathcal{F}({\rm P}_t)=\left[-\sqrt{2t+1};\sqrt{2t+1}\right]$ and $\varphi(b+td)=-\sqrt{2t+1}$. Meanwhile, $\mathcal{F}({\rm D})=\left\{y=(1,y_2,y_3)^T \in \mathbb{R}^3 \; : \;  y_3 \geq \sqrt{1+y_2^2}\right\}$ and $\mathcal{S}({\rm D})=\left\{(1,0,1)^T\right\}$. Therefore, 
		$\varphi(b)+ t \sup \{d^Ty \; : \; y \in \mathcal{F}{\rm (D)}\}=-1-t$ and $$\varphi(b)+ t \sup \{d^Ty \; : \; y \in \mathcal{S}{\rm (D)}\}=-1.$$ It is clear that $-1-t \leq -\sqrt{2t+1} < -1$ for every $t \in  (0,1)$. 
}\end{example}

\section{Properties of the Function $\psi(\cdot)$}

In this section, differentiability propertiesof the function $\psi(\cdot)$ and some increment estimates will be obtained.

\begin{theorem}\label{directional_derivative_2} Suppose that both problems {\rm (P)} and {\rm (D)} are strictly feasible. Then, one has 
	\begin{equation}\label{diff_psi}
		\psi'(c;h)=\begin{cases}
			\inf\limits_{x \in \mathcal{S}({\rm P})} x^Th & \text{ if } \quad h \in A^T(\mathbb{R}^m)\\
			-\infty & \text{ if } \quad h \notin A^T(\mathbb{R}^m).
		\end{cases}
	\end{equation}	 
\end{theorem}
{\bf Proof} By the assumed strict feasibility of  {\rm (P)} and {\rm (D)} (see Definition~\ref{strict_feassibility}), there exist $x^0 \in \mathbb{R}^n$ and $y^0 \in \mathbb{R}^m$ such that $Ax^0 >_Kb$,  $y^0>_{K^*}0$, and $A^Ty^0=c$. Since $\psi'(c;0)=0$, formula~\eqref{diff_psi} is valid for $h=0$. Fix a vector $h \in \mathbb{R}^n\setminus\{0\}$ and define the function $g: \mathbb{R} \to \overline{\mathbb{R}}$ by
\begin{equation*}
	g(t):=-\psi(c+th)=-\inf\{(c+th)^Tx \; : \; x \in \mathbb{R}^n, \, Ax \geq_{K}b\}\quad (t \in \mathbb{R}).
\end{equation*}
Note that $g(t)\neq -\infty$ for all $t\in\mathbb R$ and $\psi'(c;h)=-g'(0;1)$. Thanks to the weak duality theorem in \cite[Theorem~2.4.1 (assertion 2)]{BN_2001}, we have
$$c^Tx^0\geq v({\rm P})\geq v({\rm D})\geq b^Ty^0.$$ It follows that $v({\rm P})>-\infty$ and $v({\rm D})<+\infty$. So, by Lemma~\ref{Strong_duallity},  {\rm (P)} and {\rm (D)} have solutions, and $v({\rm D})=v({\rm P})$. Hence, the relations $g(0)=-\psi(c)=-v{\rm (P)}$ imply that the value  $g(0)$ is finite. Clearly, $g(t)\geq -(c+th)^Tx^0$ for any $t \in \mathbb{R}$.

If $h \in A^T(\mathbb{R}^m)$, then $h=A^Tv$ for some $v \in \mathbb{R}^m$. Note that $$A^T(y^0+tv)=A^Ty^0+tA^Tv=c+th.$$ Since $y^0 >_{K^*}0$, one can find $\tau>0$ such that $y^0+tv >_{K^*}0$ for every $t \in ]-\tau, \tau[$. For any $t \in ]-\tau, \tau[$, applying the weak duality theorem in \cite[Theorem~2.4.1 (assertion~2)]{BN_2001} for the conic linear problem \begin{equation}\label{P_perturbed}\min\big\{(c+th)^Tx \; : \; x \in \mathbb{R}^n,\ \, Ax \geq_{K}b\big\}\end{equation} and its dual \begin{equation}\label{D_perturbed}\max\big\{b^Ty \; : \; y \in \mathbb{R}^m, \, A^Ty=c+th,\ \, y\geq_{K^*}0\big\},\end{equation} 
one has $(c+th)^Tx \geq b^T(y^0+tv)$ for all $x \in \mathbb{R}^n$ satisfying $Ax \geq_K b$. Consequently, $\inf\big\{(c+th)^Tx \; : \; x \in \mathbb{R}^n, \, Ax \geq_{K}b\big\} \geq b^T(y^0+tv).$ So, we have $g(t) \leq - b^T(y^0+tv)$. Since $g(t)\geq -(c+th)^Tx^0$ for all $t \in \mathbb{R}$, this implies that $g(t)$ is finite for every $t \in ]-\tau, \tau[$. Now, fixing a number $t \in ]-\tau, \tau[$ and applying Lemma~\ref{Strong_duallity} for the problems~\eqref{P_perturbed} and~\eqref{D_perturbed}, one has 
\begin{multline*}
	\inf\big\{(c+th)^Tx \; : \; x \in \mathbb{R}^n, \, Ax \geq_{K}b\}\\
	=\sup\{b^Ty \; : \; y \in \mathbb{R}^m, \, A^Ty=c+th, \, y\geq_{K^*}0\big\}.
\end{multline*}   
Therefore,
\begin{equation}\label{g(t)}
	\begin{aligned}
		g(t)&=-\sup\big\{b^Ty \; : \; y \in \mathbb{R}^m, \, A^Ty=c+th, \, y\geq_{K^*}0\big\}\\
		&=\inf\big\{(-b)^Ty \; : \; y \in \mathbb{R}^m, \, A^Ty=c+th, \, y\geq_{K^*}0\big\}.
	\end{aligned}
\end{equation}
To apply a result from~\cite{An_Yen_2015}, we define a function $f:\mathbb{R} \times \mathbb{R}^m\to\mathbb R$ by setting $$f(t, y)=(-b)^Ty\quad \forall (t,y)\in \mathbb{R} \times \mathbb{R}^m.$$ Clearly, $f$ is convex and continuous on $\mathbb{R} \times \mathbb{R}^m$. Let $G: \mathbb{R} \rightrightarrows \mathbb{R}^m$ be the multifunction given by $G(t):=\big\{y \in \mathbb{R}^m \; : \; A^Ty=c+th, \, y\geq_{K^*}0 \big\}.$ Note that $G$ is a convex multifunction. Consider the convex optimization problem 
\begin{equation}\label{parametric_new}\min \{f(t,y)\; : \; y\in G(t)\},\end{equation} which depends on the parameter $t\in\mathbb R$, and observe by~\eqref{g(t)} that $g(\cdot)$ is the optimal value function of this problem.  According to a remark given in ~\cite[p.~113]{An_Yen_2015}, $g(\cdot)$ is a convex function. Let $\bar{y}$ be a solution of $({\rm D})$. Then, it is easy to verify that $\bar y$ is a solution of the parametric problem~\eqref{parametric_new} at $\bar t:=0$. Hence,  applying~\cite[Theorem~4.2]{An_Yen_2015} for~\eqref{parametric_new} yields
\begin{equation*}
	\partial g(0)=\bigcup\limits_{(t^*,y^*) \in \partial f(0, \bar y)}   \big\{t^* + D^*G( 0, \bar y)(y^*) \big\}.
\end{equation*}
Since $\partial f(t,y)=\left\{\begin{bmatrix}
	0 \\ -b
\end{bmatrix}\right\}$, this implies that
\begin{equation}\label{subdifferential_of_g_2}\partial g(0)=D^*G(0, \bar y)(-b).
\end{equation} 
By the definition of coderivative for convex multifunctions, 
\begin{equation*}
	\begin{aligned}
		D^*G(0, \bar y)(-b)&=\{t^* \in \mathbb{R} \; : \; (t^*, b) \in N((0,\bar y), \rm{gph}\, G)\}\\
		&=\big\{t^* \in \mathbb{R} \; : \; t^* t + b^T(y - \bar y) \leq 0 \ \, \forall (t,y) \in \rm{gph}\, G\big\}.
	\end{aligned}
\end{equation*} Hence,
\begin{multline}\label{sub_diff_g_2} D^*G(0, \bar y)(-b)\\=\big\{t^* \in \mathbb{R} \; : \;  b^Ty+t^* t  \leq  b^T\bar y \ \; \forall (t,y) \textrm{ s.t. } A^Ty-th=c,\, y\geq_{K^*}0 \big\}.\end{multline}
Let
\begin{equation*}
	\widetilde{b}:=\begin{bmatrix}
		b\\
		t^*
	\end{bmatrix}, \quad \widetilde{A}:=\begin{bmatrix}
		A \\ -h^T
	\end{bmatrix}, \quad \widetilde{y}:=\begin{bmatrix}
		y\\
		t
	\end{bmatrix}, \quad \widetilde{K}:=K \times \{0\}, \quad \alpha:=b^T\bar{y}.
\end{equation*}
Since $\widetilde{K}^*=K^* \times \mathbb{R}$, the system $\begin{cases}
	A^Ty-th=c\\ y\geq_{K^*} 0
\end{cases}$ is equivalent to $\begin{cases}
	\widetilde{A}^T\widetilde{y}=c\\ \widetilde{y}\geq_{\widetilde{K}^*} 0.
\end{cases}$ Hence, from~\eqref{sub_diff_g_2} one gets
\begin{equation}\label{sub_diff_g_3}
	\begin{aligned}
		D^*G(0, \bar y)(-b)=\big\{t^* \in \mathbb{R} \; : \; \widetilde{b}^T\widetilde{y} \leq \alpha \ \forall \widetilde{y}\ \, \textrm{with}\ \, \widetilde{A}^T\widetilde{y}=c,\ \widetilde{y}\geq_{\widetilde{K}^*} 0\big\}.
	\end{aligned}
\end{equation} 
Setting $\widetilde{y}^0:=\begin{bmatrix}
	y^0\\
	0
\end{bmatrix}$, we have $\widetilde{A}^T\widetilde{y}^0 =c$ and $\widetilde{y}^0>_{\widetilde{K}^*} 0$. Hence, by the Farkas-type result in Lemma~\ref{Farkas_type}(b), we can assert that the inequality $\widetilde{b}^T\widetilde{y} \leq \alpha$ holds for every $\widetilde{y} \in  \mathbb{R}^m \times \mathbb{R}$ satisfying $\widetilde{A}^T\widetilde{y}=c,\  \widetilde{y}\geq_{\widetilde{K}^*} 0$ iff  there exists $x \in \mathbb{R}^n$ satisfying $\widetilde{A}x\geq_{\widetilde{K}}\widetilde{b}$ and $c^Tx \leq \alpha$. Since $\widetilde{A}x=\begin{bmatrix}
	Ax \\ -h^Tx
\end{bmatrix}$, one sees that the inequality $\widetilde{A}x\geq_{\widetilde{K}}\widetilde{b}$ is equivalent to the system of conditions $Ax \geq_K b$ and $-h^Tx=t^*$.
Therefore, combining~\eqref{sub_diff_g_3} with the description of the feasible region of~({\rm P}), we have
\begin{equation*}
	\begin{aligned}
		D^*G(0, \bar y)(-b)&=\{t^* \in \mathbb{R} \; : \; \exists x \in \mathbb{R}^n \; \textrm{s.t.} \; \widetilde{A}x\geq_{\widetilde{K}}\widetilde{b},\; c^Tx \leq \alpha\}\\
		&=\{t^* \in \mathbb{R} \; : \; \exists  x \in \mathbb{R}^n \; \textrm{s.t.} \; Ax \geq_K b,\; -h^Tx=t^*,\; c^Tx \leq \alpha\}\\
		&=\{-h^Tx \; : \; x \in \mathbb{R}^n,\; Ax \geq_K b,\; c^Tx\leq b^T\bar y \}\\
		&=\{-x^Th \; : \; x \in \mathcal{F}({\rm P}),\;  c^Tx\leq b^T\bar y\}.
	\end{aligned}
\end{equation*}  Hence,
\begin{equation}\label{co_diff_G_2} D^*G(0, \bar y)(-b)=\big\{-x^Th \; : \; x \in \mathcal{F}({\rm P}),\;  c^Tx\leq v({\rm D})\big\}.
\end{equation}
Since $v({\rm D})=v({\rm P})$, by~\eqref{co_diff_G_2} one has
\begin{equation*}
	D^*G(0, \bar y)(-b)=\{-x^Th \; : \; x \in \mathcal{F}({\rm P}),\;  c^Tx\leq v({\rm P})\}=\{-x^Th \; : \; x \in \mathcal{S}({\rm P})\}.
\end{equation*}   
From this and \eqref{subdifferential_of_g_2} it follows that $\partial g(0)=\{-x^Th \; : \; x \in \mathcal{S}({\rm P})\}.$
Then, since $g(\cdot)$ is a proper convex function and $0\in {\rm int} ({\rm dom}\, g)$, by \cite[Theorem~23.4]{Rockafellar_1970} we have
\begin{equation}\label{three_equalities_3}
	g'(0;1)=\sup\limits_{s \in \partial g(0)}s=\sup\limits_{x \in \mathcal{S}({\rm P})} (-x^Th).
\end{equation} As $\psi'(c;h)=-g'(0;1)$,~\eqref{three_equalities_3} yields 
\begin{equation*}\label{three_equalities_4}
	\psi'(c;h)=-\sup\limits_{x \in \mathcal{S}({\rm P})} (-x^Th)=\inf\limits_{x \in \mathcal{S}({\rm P})} x^Th.
\end{equation*}
We have thus proved the first assertion in~\eqref{diff_psi}.

Now, let $h \notin A^T(\mathbb{R}^m)$. In this case, we have $g(t)=+\infty$ for every $t \in \mathbb{R}\setminus\{0\}$. Indeed, if $g(t)\in\mathbb R$\ for some  $t \in \mathbb{R}\setminus\{0\}$, then the objective function of the minimization problem~\eqref{P_perturbed}  is bounded from below. So, applying Lemma~\ref{Strong_duallity} for the problems~\eqref{P_perturbed} and \eqref{D_perturbed}, we can assert that~\eqref{D_perturbed} has a solution $y^1$. Since $A^Ty^1=c+th$ and $A^Ty^0=c$, one gets $h=A^T\left[t^{-1}(y^1-y^0)\right]$. This is impossible because $h \notin A^T(\mathbb{R}^m)$. Thus, $g(t)=+\infty$ for every $t \in \mathbb{R}\setminus\{0\}$. It follows that $g'(0;1)=+\infty$. Then we have $\psi'(c;h)=-g'(0;1)=-\infty$.
$\hfill \Box$

\begin{proposition}\label{Thm2} For every $h \in \mathbb{R}^m$ and $t>0$, one has
	\begin{equation}\label{inc_psi_lower}
		\psi(c+th)\geq \psi(c)+t \inf\limits_{x \in \mathcal{F}{\rm (P)}} h^Tx, 
	\end{equation}	
	and
	\begin{equation}\label{inc_psi_upper}\psi(c+th) \leq \psi(c)+t \inf\limits_{x \in \mathcal{S}{\rm (P)}} h^Tx.
	\end{equation}
\end{proposition}
{\bf Proof} Let $h \in \mathbb{R}^m$ and $t>0$ be given arbitrarily. Clearly,
\begin{equation*}
	\begin{aligned}
		\psi(c+th) &=\inf\big\{(c+th)^T x \; : \; x \in \mathcal{F}{\rm (P)}\big\}\\
		&\geq \inf\big\{c^T x \; : \; x \in \mathcal{F}{\rm (P)}\big\} + \inf\big\{th^Tx \; : \; x \in \mathcal{F}{\rm (P)}\big\} \\
		&=\psi(c)+t \inf\big\{h^Tx \; : \; x \in \mathcal{F}{\rm (P)}\big\}.  
	\end{aligned}
\end{equation*} So, the inequality in~\eqref{inc_psi_lower} is valid. In addition, 
\begin{equation*}
	\begin{aligned}
		\psi(c+th) &=\inf\big\{(c+th)^T x \; : \; x \in \mathcal{F}{\rm (P)}\big\}\\
		&\leq \inf\big\{(c+th)^T x \; : \; x \in \mathcal{S}{\rm (P)}\big\}\\
		&=  \psi(c)+t \inf\limits_{x \in \mathcal{S}{\rm (P)}} h^Tx. 
	\end{aligned}
\end{equation*} This means that the estimate~\eqref{inc_psi_upper} holds.
$\hfill \Box$ 

\begin{example}\label{Exam6}{\rm Consider problem ${\rm (P)}$ in the setting and notations of Example~\ref{Exam4}. Choose $h=(0,-1)^T$. For $0<t<1$, using~\eqref{exam_feasible_region} one has 
		\begin{equation*}
			\begin{aligned}
				\psi(c+th)&=\inf\Big\{x_1-tx_2 \; : \; x=(x_1,x_2)\in\mathbb R^2,\ x_1\geq \sqrt{1+x_2^2}\Big\}\\
				&=\inf\Big\{\sqrt{1+x_2^2}-tx_2 \; : \; x_2 \in \mathbb{R}\Big\}\\
				&=\sqrt{1-t^2}.
			\end{aligned}
		\end{equation*} Note that $\psi(c)=\varphi(b)=1$. Since $\mathcal{F}({\rm P})=\Big\{x \in \mathbb{R}^2 \; : \; x_1 \geq \sqrt{1+x_2^2}\Big\}$,~\eqref{inc_psi_lower} gives the trivial lower estimate $\psi(c+th)\geq -\infty$ for every $t>0$. Meanwhile, as $\mathcal{S}({\rm P})=\big\{(1,0)^T\big\}$,~\eqref{inc_psi_upper} gives the upper estimate $\psi(c+th)\leq 1$ for every $t>0$. Therefore, both estimates~\eqref{inc_psi_lower} and~\eqref{inc_psi_upper} are strict for all $t\in ]0,1[$}.\end{example}	

The next increment formula for $\psi$ is a generalization of the corresponding result stated in~\cite[p.~119]{Gauvin_2001}, where the case $K=\mathbb{R}^m_+$ was considered.

\begin{theorem}Suppose that $K$ is a polyhedral convex cone, {\rm (P)} is feasible, and {\rm (D)} is strictly feasible. Then, for every $h \in A^T(\mathbb{R}^m)$, there exists $\tau >0$ such that 
	\begin{equation}\label{inc_psi}
		\psi(c+th) = \psi(c)+t \inf\limits_{x \in \mathcal{S}{\rm (P)}} h^Tx.
	\end{equation}
	for all $t \in ]0,\tau]$. 
\end{theorem}
{\bf Proof} By our assumptions, there exists $y^0 \in \mathbb{R}^m$ such that $A^Ty^0=c$ and $y^0>_{K^*}0$. Since $h \in A^T(\mathbb{R}^m)$, there is $v \in \mathbb{R}^m$ satisfying $h=A^Tv$. Then, one has $A^T(y^0+tv)=A^Ty^0+tA^Tv=c+th.$ As $y^0 >_{K^*}0$, there exists a number $\tau'>0$ such that $y^0+tv >_{K^*}0$ for all $t \in ]-\tau', \tau'[$. For each $t \in ]-\tau', \tau'[$, the linear optimization problem $\min\big\{(c+th)^Tx \; : \; x \in \mathbb{R}^n,\ \, Ax \geq_{K}b\big\}$ 
and its dual $\max\big\{b^Ty \; : \; y \in \mathbb{R}^m, \, A^Ty=c+th,\ \, y\geq_{K^*}0\big\}$ 
are feasible. Therefore, by the well-known strong duality theorem in linear programming one has 
\begin{equation*}
	\begin{aligned}
		\psi(c+th)&=\inf\big\{(c+th)^Tx \; : \; x \in \mathbb{R}^n,\ \, Ax \geq_{K}b\big\}\\
		&=\sup\big\{b^Ty \; : \; y \in \mathbb{R}^m, \, A^Ty=c+th,\ \, y\geq_{K^*}0\big\} 
	\end{aligned}
\end{equation*}
and $\psi(c+th)$ is finite. So, 
\begin{equation}\label{equation1}
	\begin{aligned}
		\psi(c+th)&=-\inf\big\{(-b)^Ty \; : \; y \in \mathbb{R}^m, \, A^Ty=c+th,\ \, y\geq_{K^*}0\big\}\\
		&= -\inf\big\{(-b)^Ty \; : \; y \in \mathbb{R}^m, \, A^Ty\geq c+th,\\
		&\hspace*{4cm} -A^Ty\geq -(c+th), \, y\geq_{K^*}0\big\}\\
		&= -\inf\big\{(-b)^Ty \; : \; y \in \mathbb{R}^m, \, My \geq_L \widetilde{c} + t \widetilde{h}  \big\},\\
	\end{aligned}
\end{equation}
where $$M:=\begin{bmatrix}
	A^T \\ -A^T	\\	E_m
\end{bmatrix} \in \mathbb{R}^{(2n+m)\times m}, \ \, \widetilde{c}:=\begin{bmatrix}
	c \\ -c \\ 0
\end{bmatrix} \in \mathbb{R}^{2n+m}, \ \, \widetilde{h}:=\begin{bmatrix}
	h \\ -h \\ 0
\end{bmatrix} \in \mathbb{R}^{2n+m},$$ $E_m$ denotes the unit matrix in $\mathbb R^{m\times m}$, and $L:=\mathbb{R}^n_+ \times \mathbb{R}^n_+ \times K^*$. Clearly, $L$ is a polyhedral convex cone. So, applying Theorem~\ref{PCC_Case} to the problem
$\inf\big\{(-b)^Ty \; : \; y \in \mathbb{R}^m, \, My \geq_L \widetilde{c} + t \widetilde{h}  \big\},$
one finds a number $\tau \in ]0, \tau'[$ such that for any $t \in ]0, \tau[$ one has  
\begin{multline*}
	\inf\big\{(-b)^Ty \; : \; y \in \mathbb{R}^m, \, My \geq_L \widetilde{c} + t \widetilde{h}  \big\}= \inf\big\{(-b)^Ty \; : \; y \in \mathbb{R}^m, \, My \geq_L \widetilde{c} \big\}\\
	\hspace*{1cm} \ + \  t \sup \big\{\widetilde{h}^Tz \; : \; z \in \mathbb{R}^{2n+m}, \, M^Tz=-b, \, z\geq_{L^*}0, \widetilde{c}^Tz = \inf\limits_{My \geq_L \widetilde{c}}\, (-b)^Ty \big\}.
\end{multline*}
Combining this with~\eqref{equation1}, we have
\begin{multline}\label{equation2}
	\psi(c+th)= \psi(c)\\ - t \sup \big\{\widetilde{h}^Tz \; : \; z \in \mathbb{R}^{2n+m}, \, M^Tz=-b, \, z\geq_{L^*}0, \,  \widetilde{c}^Tz = -\psi(c) \big\},
\end{multline} where the equality $\psi(c)=\inf\big\{(-b)^Ty \; : \; y \in \mathbb{R}^m, \, My \geq_L \widetilde{c} \big\}$ follows from~\eqref{equation1} if one takes $t=0$. 
Substituting $z= \begin{bmatrix}
	u \\ v \\w 
\end{bmatrix}$ with $u, v \in \mathbb{R}^n$ and $w \in \mathbb{R}^m$ into~\eqref{equation2}, one gets
\begin{equation*}\label{equation3}\begin{array}{rl}
	\psi(c+th)= \psi(c) - t \sup \big\{& h^Tu-h^Tv\;\;  : \;\  u, v \in \mathbb{R}^n_+, \, w \in K,\\	& Au-Av+w=-b,\, c^Tu-c^Tv = -\psi(c) \big\}.
\end{array}
\end{equation*}
Setting $x:=v-u$, one has
\begin{equation*}
	\begin{aligned}
		\psi(c+th)&= \psi(c) \\
		&\hspace{0.3cm}- t \sup \big\{-h^Tx\; : \; x \in \mathbb{R}^n, \, w \in K, w=Ax-b,\, -c^Tx = -\psi(c) \big\}\\
		&=\psi(c) + t \inf \big\{h^Tx\; : \; x \in \mathbb{R}^n, \, w \in K, w=Ax-b,\, c^Tx = \psi(c) \big\}\\ 
		&=\psi(c) + t \inf \big\{h^Tx\; : \; x \in \mathbb{R}^n, \, Ax \geq_K b,\, c^Tx = \psi(c) \big\}\\ 
		&= \psi(c)+t \inf\limits_{x \in \mathcal{S}{\rm (P)}} h^Tx.	
	\end{aligned}	
\end{equation*}
This proves that~\eqref{inc_psi} holds. $\hfill \Box$ 

\section{Conclusions}

Based on strict feasibility conditions and the duality theory in the conic linear programming, we have shown that two optimal value functions of a conic linear program, whether either the constraint system is linearly perturbed or the objective function is linearly perturbed, have some nice continuity and differentiability properties. It turns out that if the convex cone is non-polyhedral, then the behaviors of the optimal value functions in question are more complicated than that of the corresponding optimal value functions of a linear program.  

Properties of the solution maps of parametric conic linear programs and generalizations of our results for infinite-dimensional conic linear programming deserve further investigations.

\begin{acknowledgements}
This work was supported by National Foundation for Science $\&$ Technology Development (Vietnam) and Pukyong National University (Busan, Korea). The first author was partially supported by the Hanoi National University of Education under grant number SPHN20-07: \textit{``Differential Stability in Conic Linear Programming''}. The second author was supported by the National Research Foundation of Korea Grant funded by the Korean Government (NRF-2019R1A2C1008672).
\end{acknowledgements}


\end{document}